\documentclass[10pt,twoside]{article}
\usepackage{graphicx}
\usepackage{amsmath}
\usepackage{Latex-document}
\font\mathsymb=msbm10
\def\Bbb#1{\hbox{\mathsymb #1}}

\title{\bf Teaching Linear Algebra at University\vskip 6mm}

\author{J.-L. Dorier\vspace*{-0.5cm}\thanks{IUFM de Lyon et \'equipe DDM,
Laboratoire Leibniz, 46, ave. F. Viallet, 38 031 Grenoble Cedex,
France. E-mail: Jean-Luc.Dorier@imag.fr}}
\date{\vspace{-8mm}}

\begin{document}

\maketitle

\thispagestyle{first} \setcounter{page}{875}

\begin{abstract}

\vskip 3mm

Linear algebra represents, with calculus, the two main mathematical
subjects taught in science universities. However this teaching has
always been difficult. In the last two decades, it became an active
area for research works in mathematics education in several countries.
Our goal is to give a synthetic overview of the main results of these
works focusing on the most recent developments.
The main issues we will address concern:
\begin{itemize}
\item the epistemological specificity of linear algebra and the
interaction with research in history of mathematics
\item the cognitive
flexibility at stake in learning linear algebra
\item three principles for
the teaching of linear algebra as postulated by G. Harel
\item the
relation between geometry and linear algebra
\item an original teaching
design experimented by M. Rogalski
\end{itemize}

\vskip 4.5mm


\noindent {\bf Keywords and Phrases:} University teaching, Linear
algebra, Curriculum,
  Vector space, Representation, Geometry, Cognitive flexibility,
Epistemology.
\end{abstract}

\vskip 12mm

\section{Introduction} \label{section 1}\setzero
\vskip-5mm \hspace{5mm}

In most countries, science-orientated curricula in the first
two years at university consist of courses in two main subjects, namely,
calculus and linear algebra. The difficulties in these two fields are of
different nature. Mathematics education research first developed works on
calculus, but in the past 20 years, an increasing number of studies has
been carried out about the teaching of linear algebra. One can
distinguish roughly two main traditions in the teaching of linear
algebra: one focuses on the study of formal vector spaces while the other
proposes a more analytical approach based on the study of $\Bbb R^n$
and matrix calculus. Between these two orientations, there exist a
continuum of teaching designs, in which each pole is more or less
dominant. However, the teaching of linear algebra is universally
recognised as difficult. Students usually feel that they land on another
planet, they are overwhelmed by the number of new definitions and the
lack of connection with previous knowledge. On the other hand, teachers
often feel frustrated and disarmed when faced with the inability of their
students to cope with ideas that they consider to be so simple. Usually,
they incriminate the lack of practice in basic logic and set theory or the
impossibility for the students to use geometrical intuition. These
complaints have a certain validity, but the few attempts at remedying
this state of affairs - with the teaching of Cartesian geometry or/and
logic and set theory prior to the linear algebra course - did not seem to
improve the situation substantially. The aim of this text is to give an
account of the main trends in this area of mathematics education research.

\section{Historical analyses} \label{section 2}
\setzero\vskip-5mm \hspace{5mm }

An epistemological analysis of the history of linear algebra is a way to reveal some possible sources of students'
difficulties as well an inspiration in the design of activities for students. Several works have been carried out
in this direction (see [4], [6], [8], [17] and [24]). In this paper, we will give an account of only one of the
main result of this type of research. It concerns the last phase of the genesis of the theory of vector spaces,
whose roots can be found in the late nineteenth century, but really started only after 1930. It corresponds to the
axiomatisation of linear algebra, that is to say a theoretical reconstruction of the methods of solving linear
problems, using the concepts and tools of a new axiomatic central theory. These methods were operational but they
were not explicitly theorised or unified. It is important to realise that this axiomatisation did not, in itself,
allow mathematicians to solve new problems; rather, it gave them a more universal approach and language to be used
in a variety of contexts (functional analysis, quadratic forms, arithmetic, geometry, etc.). The axiomatic
approach was not an absolute necessity, except for problems in non-denumerable infinite dimension, but it became a
universal way of thinking and organising linear algebra. Therefore, the success of axiomatisation did not come
from the possibility of reaching a solution to unsolved mathematical problems, but from its power of
generalisation and unification and, consequently, of simplification in the search for methods for solving problems
in mathematics.

As a consequence, one of the most noticeable difficulties
encountered in the learning of unifying and generalising concepts
are associated with the pre-existing, related elements of
knowledge or competencies of lower level.  Indeed, these need to
be integrated within a process of abstraction, which means that
they have to be looked at critically, and their common
characteristics have to be identified, and then generalised and
unified. From a didactic point of view, the difficulty is that any
linear problem within the reach of a first year university student
can be solved without using the axiomatic theory. The gain in
terms of unification, generalisation and simplification brought by
the use of the formal theory is only visible to
the expert.

  One solution would be to give up teaching the formal theory
of vector spaces. However, many people find it important that students
starting university mathematics and science studies get some idea about
the axiomatic algebraic structures of which vector space is one of the
most fundamental. In order to reach this goal, the question of formalism
cannot be avoided. Therefore, students have to be introduced to a certain
type of reflection on the use of their previous elements of knowledge and
competencies in relation with new formal concepts. This led Dorier,
Robert, Robinet and Rogalski to introduce what they called {\it `meta
level activities'} (see [5], [6], [9], [11], [21] [22] and [23]). These
activities are introduced and maintained by an explicit discourse on the
part of the teacher about the significance of the introduced concepts for
the general theory, their generalising and unifying character, the change
of point of view or a theoretical detour that they offer, the types of
general methods they lead to, etc. It hinges on the general attitude of
the teacher who induces a constant underlying meta-questioning concerning
new possibilities or conceptual gains provided by the use of linear
algebra concepts, tools and methods.

\section{Cognitive flexibility} \label{section 3} \setzero\vskip-5mm
\hspace{5mm }

One of the main difficulties in learning linear algebra has to do
  with the variety of languages, semiotic registers of representation,
points of view and settings through which the objects of linear
algebra can be represented. Students have to distinguish these
various ways of representing objects of linear algebra, but they
also need to translate from one to another type and, yet, not
confuse the objects with their different representations. These
abilities could be referred by the general notion of cognitive
flexibility. This question is central in several works on the
teaching and learning of linear algebra.

Students' difficulties
with the formal aspect of the theory of vector spaces are not just
a general problem with formalism but mostly a difficulty of
understanding the specific use of formalism within the theory of
vector spaces and the interpretation of the formal concepts in
relation with more intuitive contexts like geometry or systems of
linear equations, in which they historically emerged. Various
diagnostic studies conducted by Dorier, Robert, Robinet and
Rogalski pointed to a single massive obstacle appearing for all
successive generations of students and for nearly all modes of
teaching, namely, what these authors termed {\it the obstacle of
formalism} (see [6], [7], [10] and  [27]).

In [16], Hillel distinguished three basic languages used in linear algebra:
  the `abstract language' of the general abstract theory, the `algebraic
language' of the $\Bbb R^n$ theory and the  `geometric language'
of the two- and three-dimensional spaces. The `opaqueness' of the
representations seems to be ignored by lecturers, who constantly
shift the notations and modes of description, without alerting the
students in any explicit way. By far, the most confusing case for
students is the shift from the abstract to the algebraic
representation when the underlying vector space is $\Bbb R^n$. In
this case, an $n$-tuple (or a matrix) is represented as another
$n$-tuple (or matrix) relative to another basis. This confusion
leads to persistent mistakes in students' solutions related to
reading the values of a linear transformation given by a matrix in
a basis (see [15]). Parallel with the three languages identified
by Hillel, Sierpinska et al. distinguish three modes of thinking
that have led to the development of these languages and are
necessary for an understanding of the domain: the
`synthetic-geometric', `analytic-arithmetic', and
`analytic-structural' (see [29]).

In [12], Duval defined {\it semiotic representations} as productions
  made by the use of signs belonging to a system of representation which
has its own constraints of meaning and functioning. Semiotic
representations are, according to him, absolutely necessary in
mathematical activity, because its objects cannot be directly
perceived and must, therefore, be represented. Moreover, semiotic
representations play an essential role in developing mental
representations, in accomplishing different cognitive functions
(objectification, calculation, etc.), as well as in producing
knowledge.

In her work, Pavlopoulou (see [8] pp. 247-252) applied and tested
Duval's theory in the context of linear algebra. She distinguished
between three registers of semiotic representation of vectors: the
graphical register (arrows), the table register (columns of
coordinates), and the symbolic register (axiomatic theory of
vector spaces). Through several studies, she has shown that the
question of registers, especially as regards conversion, is not
usually taken into account either in teaching or in textbooks. She
also identified a number of student's mistakes that could be
interpreted as a confusion between an object and its
representation (especially a vector and its geometrical
representation) or as a difficulty in converting from one register
to another.

The research of Alves-Dias (see [1] and [8] pp. 252-256), an
extension of Pavlopoulou's, generalised the necessity of
conversions from one semiotic register to another for the
understanding of linear algebra to the necessity of `cognitive
flexibility'. Moreover, on the basis of Rogalski's previous work
(see [25], [26], [27] and [28]), she focused her study on the
question of articulation between the Cartesian and parametric
representations of vector subspaces, which is not a mere question
of change of register, but deals with more complex cognitive
processes involving the use of concepts like rank and duality.
Indeed, when a subspace $V$ is represented by Cartesian equations,
finding a parametric representation of $V$ mostly consists in
finding a set of generators of $V$, which is not just a change of
register, nor an elementary cognitive process, even if it is much
easier when the dimension `$d$' of $V$ is known. In any case,
competencies with regard to the concept of rank and duality are
indispensable. Moreover, in order to avoid easy mistakes in
calculations or reasoning, it is necessary to be able to have some
control over the results obtained.  Alves Dias showed that in
textbooks and classes, in general, the tasks offered to students
are very limited in terms of flexibility. She developed a series
of exercises that required the student to mobilise more changes of
settings or registers and to exert explicit control via the
concepts of rank and duality. Her experimentations demonstrated a
variety of difficulties for the students. For instance, students
often identified one type of representation exclusively through
semiotic characteristics (a representation with $x$'s and $y$'s
would be considered as obviously Cartesian) without questioning
the meaning of the representation. Concerning the means of control
over the validity of the statements by the students and
anticipation of results or answers to problems, she found that a
theorem like: $dim E = dim Ker f + dim Im f$, is known and used
correctly by many students, but it is very seldom used for those
purposes even in cases in which it would immediately bring up a
contradiction with the result obtained, or in cases in which it
offered valuable information in order to anticipate the correct
answer.

In [15], Hillel and Sierpinska stressed that a linear algebra
course which is theoretically rather than computationally framed requires
a level of thinking that is based on what has been termed by Piaget and
Garcia as the `trans-object level of analysis' which consists in the
building of conceptual structures out of what, at previous levels, were
individual objects, actions on these objects, and transformations of both
the objects and actions (see [18] p. 28). A similar claim was made by
Harel in [14], in his assertions that a substantial range of mental
processes must be encapsulated into conceptual objects by the time
students get to study linear algebra. The difficulty of thinking at the
trans-object level leads some students to develop `defense mechanisms'
(to `survive' the course), consisting in trying to produce a written
discourse formally similar to that of the textbook or of the lecture but
without grasping the meaning of the symbols and the terminology. This
appeared as a major problem for Sierpinska, Dreyfus and Hillel, and the
team set out to design an entry into linear algebra that would make this
behaviour or attitude less likely to appear in students (see [30]). The
designed teaching-learning situations were set in a dynamic geometry
environment (Cabri-geometry II) extended by several macro-constructions
for the purposes of representing a two-dimensional vector space and its
transformations (see [31] and [32]). Further analysis of the students'
behaviour in the experimented situations led Sierpinska to postulate
certain features of their thinking that could be held partly responsible
for their erroneous understandings and difficulties in dealing with
certain problems (especially the problem of extending a transformation of
a basis to a linear transformation of the whole plane). She proposed that
these features be termed ``a tendency to think in `practical' rather than
`theoretical' ways" (see [32]). The distinction between these two ways of
thinking was inspired by the Vygotskian notion of scientific, as opposed
to spontaneous or everyday concepts. The behaviour of students who were
encountering difficulties in the experimentations suggested that their
ways of thinking had the features of practical thinking rather than
theoretical thinking. In particular they had trouble going beyond the
appearance of the graphical and dynamic representations in Cabri that
they were observing and manipulating: their relation to these
representations was `phenomenological' rather than `analytic'. By far the
most blatant feature of the students' practical thinking was their
tendency to base their understanding of an abstract concept on
`prototypical examples' rather than on its definition. For example,
linear transformations were understood as `rotations, dilations, shears
and combinations of these'. This way of understanding made it very
difficult for them to see how a linear transformation could be determined
by its value on a basis, and consequently, their notion of the matrix of
a linear transformation remained at the level of procedure only.

\section{Three principles for the teaching of linear algebra}
\label{section 4} \setzero\vskip-5mm \hspace{5mm }

In [14], Harel posits three `principles' for the teaching of
linear algebra, inspired by Piaget's psychological theory of
concept development:  the {\it Concreteness Principle}, the {\it
Necessity Principle} and the {\it Generalisability Principle}.

The Concreteness Principle states, ``For students to abstract a
mathematical structure from a given model of that structure, the
elements of that model must be conceptual entities in the
student's eyes; that is to say, the student has mental procedures
that can take these objects as inputs". This principle is violated
whenever the general concept of vector space is taught as a
generalisation from less abstract structures, to students who have
not (yet) constructed the elements of these structures as mental
entities on which other mental operations can be performed.
Starting from the premise that students build their understanding
of a concept in a context that is concrete to them, Harel conclude
that a sustained emphasis on a geometric embodiment of abstract
linear algebra concepts produce a quite solid basis for students'
understanding. He insisted, however, that it would be incorrect to
conclude that a linear algebra course should start with geometry
and build the algebraic concepts through some kind of
generalisation from geometry. A teaching experiment built on this
premise allowed Harel to observe that when geometry is introduced
before the algebraic concepts have been formed, many students
remain in the restricted world of geometric vectors, and do not
move up to the general case.

The Necessity Principle --- For students to learn, they must see
an (intellectual, as opposed to social or economic) need for what
they are intended to be taught --- is based on the Piagetian
assumption (which has also been adopted by the Theory of Didactic
Situations elaborated by Brousseau in [2]) that knowledge develops
as a solution to a problem. If the teacher solves the problems for
the students and only asks them to reproduce the solutions, they
will learn how to reproduce teacher's solutions, not how to solve
problems. Deriving the definition of vector space from a
presentation of the properties of $\Bbb R^n$ is an example of a
violation of the necessity principle.

The last, Generalisability Principle postulated by Harel, is
concerned more with didactic decisions regarding the choice of
teaching material than with the process of learning itself. ``When
instruction is concerned with a `concrete' model, that is a model
that satisfies the Concreteness Principle, the instructional
activities within this model should allow and encourage the
generalisability of concepts." This principle would be violated if
the models used for the sake of concretisation were so specific as
to have little in common with the general concepts they were aimed
at. For example, the notion of linear dependence introduced in a
geometric context defined through collinearity or co-planarity is
not easily generalisable to abstract vector spaces. Harel's work
inspired curriculum reform in the US (see [3]), as well as
textbook authors (see [33]).

\section{Geometry and linear algebra}
\label{section 5} \setzero\vskip-5mm
\hspace{5mm }

In [19], Robert, Robinet and Tenaud designed and experimented
with a geometric entry into linear algebra. The aim was to overcome the
obstacle of formalism by giving a more `concrete' meaning to linear
algebra concepts, in particular, through geometrical figures that could
be used as metaphors for general linear situations in more elaborate
vector spaces. However, as Harel noticed after them in his study
mentioned above, the connection with geometry proved to be problematic.
Firstly, geometry is limited to three dimensions and therefore some
concepts, like rank, for instance, or even linear dependence, have a
quite limited field of representation in the geometric context. Moreover,
it is not rare that students refer to affine subspaces instead of vector
subspaces when working on geometrical examples within linear algebra.

  In
her work, Gueudet-Chartier (see [13] and [8] pp. 262-264) conducted an
epistemological study of the connection between geometry and linear
algebra, using the evidence from both historical and modern texts. She
found that the necessity of geometric intuition was very often postulated
by textbooks or teachers of linear algebra. However, in reality, the use
of geometry was most often very superficial. Moreover, some students
would use geometrical representations or references in linear algebra,
without this always being to their advantage. Indeed, some of them could
not distinguish the affine space from the vector space structure; they
also often could not imagine a linear transformation that would not be a
geometric transformation. In other words, the geometrical reference acted
as an obstacle to the understanding of general linear algebra. On the
other hand, some very good students were found to use geometric
references very rarely. They could operate on the formal level without
using geometrical representations. It seems that the use of geometrical
representations or language is very likely to be a positive factor, but
it has to be controlled and used in a context where the connection is
made explicit.

\section{An original teaching experiment}
\label{section 6} \setzero\vskip-5mm
\hspace{5mm }

Most of the research conducted in France on the teaching and
learning of linear algebra has been more or less directly
connected with an experimental course implemented by Rogalski (see
[9], [25], [26], [27] and [28]). This course was built on several
interwoven and long-term strategies, using meta level activities
as well as changes of settings (including intra-mathematical
changes of settings), changes of registers and points of view, in
order to obtain a substantial improvement in a sufficient number
of students.  `Long-term strategy' (see [20]) refers to a type of
teaching that cannot be divided into separate and independent
modules. The long-term aspect is vital because the mathematical
preparation and the changes in the `didactic contract' (see [2])
have to operate over a period which is long enough to be efficient
for the students, in particular as regards assessment. Moreover,
the long-term strategy refers to the necessity of taking into
account the non-linearity of the teaching due to the use of change
in points of view, implying that
a subject is (re)visited several times in the course of the year.

Rogalski's teaching design has the following main characteristics:
\begin{itemize}
  \item In
order to take into account the specific epistemological nature of the
concepts, some activities are introduced, at a favourable and precise
time of the teaching, in order to induce a reflection on a `meta' level.
\item A fairly long preliminary phase precedes the actual teaching of
elementary concepts of linear algebra. It prepares the students to
understand, through `meta' activities, the unifying role of these
concepts.
\item As much as possible, changes of settings and points of view
are used explicitly and are discussed.
\item Finally, the concept of rank is
given a central position in this teaching.
\end{itemize}
  For a long-term teaching
design, it is difficult to choose the time suitable for its evaluation,
as interference may occur due to students' own organisation of their time
and work, in a way that cannot be kept under control. Thus, phenomena of
maturing, depending on students' level of involvement (which varies
during the year) are difficult to take into account in the evaluation of
the teaching. Moreover, such a global teaching design cannot be evaluated
by usual comparative analyses, because the differences with the standard
course are too important. However, internal evaluations have been
conducted, showing several positive effects, even if some questions
remain open.

\section{Conclusions}
\label{section 7} \setzero\vskip-5mm
\hspace{5mm }

Mathematics education research cannot give a miraculous solution
to overcome all the difficulties in learning and teaching linear algebra.
Various works have consisted in diagnoses of students' difficulties,
epistemological analyses and experimental teaching, offering local
remediation. Nevertheless, these works lead to new questions, problems
and difficulties. Yet, this should not be interpreted as a failure.
Improving the teaching and learning of mathematics cannot consist in one
remediation valid for all. cognitive processes and mathematics are far
too complex for such an idealistic simplistic view. It is a deeper
knowledge of the nature of the concepts, and the cognitive difficulties
they enclose, that helps teachers make their teaching richer and more
expert; not in a rigid and dogmatic way, but with flexibility. In this
sense, in several countries, mathematics education research has
influenced curriculum reforms, in an non-formal way, or sometimes very
officially, like in North-America through the Linear Algebra Curriculum
Study Group (see [3]).

\label{lastpage}

\end{document}